\documentstyle{amsppt}
\pagewidth{5.30in} \pageheight{8.0in} \NoRunningHeads
\magnification = 1200

\topmatter
\title Rigidity Theorems For Lagrangian Submanifolds of $C^n$ and $CP^n$
With Conformal Maslov Form\endtitle
\author Xiaoli Chao and Yuxin Dong*\endauthor
\thanks {*Supported by Zhongdian grant
of NSFC}
\endthanks
\abstract {In this paper, we obtain a rigidity theorem for
Lagrangian submanifolds of  $C^n$ and $CP^n$  with conformal
Maslov form.}
\endabstract
, {53C42}
\endsubjclass
\endtopmatter
\document
\heading{\bf 1. Introduction}
\endheading
\vskip 0.3 true cm

Let $(\widetilde{M}^{n}, J, g)$ be a K$\ddot{a}$hler manifold of
complex dimension $n$. The K$\ddot{a}$hler form $\omega$ on
$\widetilde{M}^{n}$ is given by $\omega (X,Y)=g( X, JY )$.  An
immersion $\psi:M^n\rightarrow \widetilde{M}^{n}$ of an
$n$-dimensional manifold $M$ is called Lagrangian if
$\psi^{\ast}\omega \equiv 0$. Lagrangian submanifolds in a
K$\ddot{a}$hler manifold or more generally in a symplectic
manifold appear naturally in the context of mathematical physics.
Since 1970s, various kind of Lagrangian submanifolds in
K$\ddot{a}$hler manifolds have also been investigated extensively
from Riemannian geometric point of view (see [1] and the
references therein).

From [4], we know that the $n$-sphere cannot be embedded in $C^n$
as a Lagrangian submanifold. This result is not true when the
Lagrangian sphere is immersed but not embedded. The simplest
immersions of $S^n$ into $C^n$, known as Whitney spheres [8], are
induced by a map $\psi: E^{n+1}\rightarrow C^n (\cong R^{2n})$
defined by
$$
\psi(x_0,x_1,\cdots ,x_n)=\frac{r}{1+x^2_0}
(x_1,\cdots,x_n,x_0x_1,\cdots ,x_0x_n)+A
$$
where $r$ is a positive number and $A$ is a vector in $C^{n}$. We
will refer to $r$ and $A$ as the radius and the center of the
Whitney sphere respectively.  The Whitney spheres have the best
possible behavior either from the viewpoint of topology or the
viewpoint of submanifold geometry. They have a unique
self-intersection point $\psi(-1,0,\cdots,0)=\psi(1,0,\cdots,0)$.
Their second fundamental form $h$ satisfy
$$
h(X,Y)=\frac{n}{n+2}\{\langle X,Y\rangle H+\langle JX,H\rangle
JY+\langle JY,H\rangle JX\}\tag{1}
$$
where $X,Y$ are tangent to $\psi$ and $H=\frac{1}{n}trace(h)$ is
the mean curvature vector of the Whitney sphere. This property may
be regarded as the Lagrangian version of umbilicity. In [7], it
was proved that the Whitney spheres are the only closed Lagrangian
submanifolds in $C^n$ having this property.

From [1], we know that there is no Lagrangian immersion of $S^n$
into $C^n$ with parallel mean curvature vector. However, the
Whitney spheres have the property that $JH$ are conformal vector
fields. We may regard a Lagrangian submanifold with this property
as the analogue of hypersurfaces of constant mean curvature. As
the dual form of $JH$ is the Maslov form of the Lagrangian
immersion, the Lagrangian submanifolds whose $JH$ are conformal
vector field will be known as Lagrangian submanifolds with
conformal Maslov form (see [7]). In [7], these submanifolds were
studied when the ambient space is $C^n$. The authors proved in [7]
that the Whitney spheres are the only compact Lagrangian
submanifolds of $C^n$ with conformal Maslov form and the null
first Betti number. This result may be regarded as the Lagrangian
version of the classical Hopf theorem. When the ambient spaces are
$CP^n$ and $CH^n$, the authors in [2], [4] constructed similar
Lagrangian immersions from $S^n$ into these spaces, which were
called Whitney spheres too. It turns out that (1) also
characterizes Whitney spheres among closed non-minimal Lagrangian
submanifolds when the ambient spaces are $C^n$ and $CP^n$. In the
case of $CH^n$, there are other two families of examples, besides
the Whitney spheres of $CH^n$, that satisfy (1).

In this paper, we first introduce a modified second fundamental
form $B$ suggested by (1). In terms of the square norm of $B$, we
will establish rigidity theorems for Lagrangian submanifolds of
$C^n$ and $CP^n$ with conformal Maslov form, and thus characterize
the Whitney spheres (Theorem 3.3) in these spaces.

\heading{\bf 2. Preliminaries}
\endheading
\vskip 0.3 true cm

Let $\psi:M^n \rightarrow \widetilde{M}^n$ be a Lagrangian
immersion of an $n$-dimensional manifold $M^n$ into a
K$\ddot{a}$hler manifold $\widetilde{M}^n$. Let $\nabla$,
$\widetilde{\nabla}$ denote the Levi-Civita connections of $M$ and
$\widetilde{M}$ respectively. The second fundamental form $h$ of
$\psi$ is given by
$$
h(X,Y)=\widetilde{\nabla}_X Y-\nabla_X Y
$$
for any vector fields $X$, $Y$ tangent to $M$. The mean curvature
vector $H$ of $\psi$ is defined by
$$
H=\frac{1}{n}trace(h)=\sum_{i=1}^n h(e_i,e_i)
$$
where $\{e_1,\cdots ,e_n\}$ is an orthonormal tangent frame of
$M$.

We choose a local orthonormal frame field $\{e_1,\cdots,e_n\}$ of
$M$ and hence $\{e_{1^*}=Je_1,\cdots,e_{n^*}=Je_n\}$ forms a local
normal frame field of $M$ in $\widetilde{M}^n$.  Let
$$\{\omega^1,\cdots,\omega^n, \omega^{1^*},\cdots,
\omega^{n^*}\}$$ denote the dual frame field of
$\{e_1,\cdots,e_n\,e_{1^*},\cdots,e_{n^*}\}$.

In this paper, we shall make use of the following convention of
indices:$$1\leq A,B,C,\cdots \leq 2n, \quad  1\leq i,j,k,m,l\cdots
\leq n$$ The structure equations of $\widetilde{M}$ are given by
$$
\aligned &d\omega^A=-\sum_{B}\omega^A_B\wedge \omega^B,\quad
\omega^A_B+\omega^B_A=0\\
&d\omega^A_B=-\sum_{C}\omega^A_C\wedge\omega^C_B+\Omega^A_B,\quad
\Omega^A_B=\sum_{CD}\frac{1}{2}K^A_{BCD}\omega^C\wedge\omega^D\\
&K^A_{BCD}+K^A_{BDC}=0\endaligned
$$
Restricting these forms to $M$, we have
$$
\aligned &\omega^i_j+\omega^j_i=0,\quad \omega^i_j=\omega^{
i^*}_{j^*},\quad \omega^{i^*}_j=\omega^{j^*}_i ,\\
&\omega^{m^*}_i=h^{m*}_{ij}\omega^j,\quad h^{m^*}_{ij}=h^{m^*}_{ji}=h^{i^*}_{mj}=h^{j^*}_{im},\\
&d\omega^i=-\omega^i_j\wedge \omega^j,\quad\omega^i_j+\omega^j_i=0\\
&d\omega^i_j=-\omega^i_k\wedge\omega^k_j+\frac{1}{2}R^i_{jkl}\omega^k\wedge\omega^l,\quad
R^i_{jkl}+R^i_{jlk}=0\\
&R^i_{jkl}=K^i_{jkl}+\sum_{m}(h_{jl}^{m^*}h_{ik}^{m^*}-h_{il}^{m^*}h_{jk}^{m^*})\\
&R^{i^*}_{j^*kl}=K^{i^*}_{j^*kl}+\sum_{m}(h_{lm}^{j^*}h_{mk}^{i^*}-h_{ml}^{i^*}h_{mk}^{j^*})
\endaligned
$$
Define the first covariant derivative of $h_{ij}^{m^*}$ by
$$
\sum_k h_{ijk}^{m^*}\omega^k=dh_{ij}^{m^*}-\sum_{k}
h_{kj}^{m^*}\omega_i^k-\sum_{k} h_{ik}^{m^*}\omega_j^k+\sum_{k}
h_{ij}^{k^*}\omega_{k^*}^{m^*}
$$
We denote by $\widetilde{M}^n(4c)$ the complex space form of
constant holomorphic sectional curvature $4c$, i.e., its curvature
satisfies
$$K^{A}_{BCD}=c\{\delta_{AC}\delta_{BD}-\delta_{AD}\delta_{BC}+J_{AC}J_{BD}-J_{AD}J_{BC}+2J_{AB}J_{CD}\}$$
where
$$(J_{AB})=\pmatrix
0 & -I_n \\
I_n & 0 \\
\endpmatrix$$
It is known that
$$
\widetilde{M}^n(4c)\cong\cases C^n, \quad &\text{if $c=0$}\\
CP^n, \quad &\text{if $c>0$}\\
CH^n, \quad &\text{if $c<0$}
\endcases
$$

For a Lagrangian submanifold $M$ in $\widetilde{M}^n(4c)$, we
obtain the Codazzi equation by a direct computation
$$h_{ijk}^{m^*}=h_{ikj}^{m^*}$$
Notice that $JH$ is a tangent vector field when the submanifold is
Lagrangian. \proclaim{Lemma 2.1}([7]) $JH$ is a conformal vector
field if and only if
$$
\nabla JH=div(JH/n)I ,
$$
which is equivalent to
$$
\sum_k h_{kk,l}^{m^*}=-div(JH)\delta_l^m
$$
\endproclaim

\definition{Definition 2.2} A Lagrangian submanifold $M$ is said to have conformal Maslov form if $JH$ is a conformal vector field of
$M$.
\enddefinition

In the introduction, we present the Whitney spheres in $C^n$. From
[4], we obtain the following family of Lagrangian immersions
$$
\psi _\theta :S^n\rightarrow CP^n,\quad\theta >0
$$
given by
$$
\psi _\theta
(x_1,...,x_n,x_{n+1})=[(\frac{(x_1,...,x_n)}{ch_\theta +ish_\theta
x_{n+1}},\frac{sh_\theta ch_\theta
(1+x_{n+1}^2)+ix_{n+1}}{ch_\theta ^2+sh_\theta ^2x_{n+1}^2})]
$$
that are called the Whitney spheres in $CP^n$. There are Whitney
immersions of $S^n$ into $CH^n$ too (see [4] for details).

\proclaim{Lemma 2.3}([1],[4], [7]) Let $\psi:M\rightarrow
\widetilde{M}^n(4c)$ be a Lagrangian immersion of an
$n$-dimensional manifold $M$ in the complex space form
$\widetilde{M}^n(4c)$ with $c\geq 0$. Then the second fundamental
form $h$ of $\psi$ satisfies
$$
h(X,Y)=\frac{n}{n+2}\{\langle X,Y\rangle H+\langle JX,H\rangle
JY+\langle JY,H\rangle JX\}
$$ for any vectors $X$ and $Y$ tangent
to $M$ if and only if either $M$ is a totally geodesic submanifold
or it is an open set of the Whitney Sphere.
\endproclaim
\remark{Remark 2.1}In the case $c<0$, there are other two families
of examples, besides Whitney spheres of $CH^n$, whose second
fundamental forms satisfy the property in Lemma 2.3. \endremark

\heading{\bf 3. Rigidity theorems for Lagrangian submanifolds}
\endheading
\vskip 0.3 true cm

In this section, we will establish a rigidity theorem for
Lagrangian submanifolds in $C^n$ and $CP^n$ with conformal Maslov
form.

Let $\psi:M^n \rightarrow \widetilde{M}^n(4c)$ be a Lagrangian
submanifold with conformal Maslov form. We define
$$
B(X,Y)=h(X,Y)-\frac{n}{n+2}\{\langle X,Y\rangle H+\langle
JX,H\rangle JY+\langle JY,H\rangle JX\}\tag{2}
$$
for any tangent vector fields $X$, $Y$ of $M$. It is easy to
verify that $trace(B)=0$. With respect to the above orthonormal
frame field $\{e_1,\cdots,e_n,e_{1^*},\cdots,e_{n^*}\}$ along $M$,
we have
$$
\aligned h_{ij}^{m^*} &=
b_{ij}^{m*}+\frac{n}{n+2}\{H^{m^*}\delta_{ij}+\langle H,
Je_i\rangle\delta_{jm}
+\langle H, Je_j\rangle\delta_{im}\}\\
&=
b_{ij}^{m*}+\frac{n}{n+2}\{H^{m^*}\delta_{ij}+H^{i^*}\delta_{jm}+H^{j^*}\delta_{im}\}\\
&= b_{ij}^{m^*}+c_{ij}^{m^*}\endaligned
$$
where
$$
c_{ij}^{m^*}=\frac{n}{n+2}\{H^{m^*}\delta_{ij}+H^{i^*}\delta_{jm}+H^{j^*}\delta_{im}\}
$$
By a direct computation, we may get
$$
|B\|^2=\|h\|^2-\frac{3n^2}{n+2}|H|^2\tag{3}
$$
Define the first and the second covariant derivatives of
$b_{ij}^{m^*}$ as follows:
$$
\aligned
{\sum_{k} b_{ijk}^{m^*}\omega^k}&=db_{ij}^{m^*}-\sum_{k}
b_{kj}^{m^*}\omega_i^k-\sum_{k} b_{ik}^{m^*}\omega_j^k+\sum_{k}
b_{ij}^{k^*}\omega_{k^*}^{m^*} \\
{\sum_l b_{ijkl}^{m^*}\omega^l}&=db_{ijk}^{m^*}-\sum_{l}
b_{ljk}^{m^*}\omega_i^l-\sum_{l} b_{ilk}^{m^*}\omega_j^l-\sum_{l}
b_{ijl}^{m^*}\omega_k^l+\sum_{l} b_{ijk}^{l^*}\omega_{l^*}^{m^*}
\endaligned
$$
\proclaim{Lemma 3.1} The tensor field $B$ is a Codazzi tensor
field, i.e., $ b_{ijk}^{m*}=b_{ikj}^{m*}$.
\endproclaim
\demo{Proof}For any point of $M$, we may take a local orthonormal
frame $\{e_i\}_{i=1}^n$ near $p$ such that $(\nabla e_i)_p=0$ for
$i=1,...,n$, which imply that $(\nabla^{\bot} e^*_i)_p = 0$ for
$i=1,...,n$. By using Lemma 2.1, we have
$$
\aligned
b_{ijk}^{m^*}&=
h_{ijk}^{m^*}-\frac{n}{n+2}\{H^{m^*}_{,k}\delta_{ij}+\langle
H_{,k}, Je_i\rangle\langle Je_j, e_m^{*}\rangle
+\langle H_{,k}, Je_j\rangle\langle Je_i, e_m^{*}\rangle\}\\
&=h_{ijk}^{m^*}-\frac{n}{n+2}\{H^{m^*}_{,k}\delta_{ij}+H^{i^*}_{,k}\delta_{jm}+H^{j^*}_{,k}\delta_{im}\}\\
&=h_{ijk}^{m^*}+\frac{div(J H)}{n+2}\{\delta_{km}\delta_{ij}+
\delta_{ik}\delta_{jm}+\delta_{jk}\delta_{im}\}\endaligned
$$
Therefore $b_{ijk}^{m^*}=b_{ikj}^{m^*}$. \qed
\enddemo

The Laplacian of $b_{ij}^{m^*}$ is defined by
$$
\triangle
b_{ij}^{m^*}=\sum_{k}b_{ijkk}^{m^*}=\sum_{k}b_{ikjk}^{m^*}=\sum_{k}b_{kijk}^{m^*}\tag{4}
$$
By Ricci identities and Lemma 3.1, we obtain
$$
\aligned
b_{kijk}^{m^*}&=b_{kikj}^{m^*}+\sum_{t}b_{kt}^{m^*}R_{ijk}^{t}+\sum_{t}b_{ti}^{m^*}R_{kjk}^{t}-
\sum_{t}b_{ki}^{t^*}R_{t^*jk}^{m^*} \\
&=b_{kki
j}^{m^*}+\sum_{t}b_{kt}^{m^*}R_{ijk}^{t}+\sum_{t}b_{ti}^{m^*}R_{kjk}^{t}-
\sum_{t}b_{ki}^{t^*}R_{t^*jk}^{m^*}\endaligned\tag{5}
$$
Now we compute
$$\frac{1}{2}\triangle\|B\|^2
=\sum_{ijkm}(b_{ijk}^{m^*})^2+\sum_{ijm} b_{ij}^{m^*}\triangle
b_{ij}^{m^*}$$ From (4),(5), we have
$$
\sum_{ijm} b_{ij}^{m^*}\triangle
b_{ij}^{m^*}=\underbrace{\sum_{ijkmt}
b_{ij}^{m^*}b_{kt}^{m^*}R_{ijk}^{t}}_{(I)}+\underbrace{\sum_{ijkmt}
b_{ij}^{m^*}b_{ti}^{m^*}R_{kjk}^{t}}_{(II)}-
\underbrace{\sum_{ijkmt}
b_{ij}^{m^*}b_{ki}^{t^*}R_{t^*jk}^{m^*}}_{(III)}
$$
where
$$
\aligned
(I)&=c\sum_{ijkmt}
b_{ij}^{m^*}b_{kt}^{m^*}(\delta_{tj}\delta_{ik}-\delta_{tk}\delta_{ij})
+\sum_{ijkmtl} b_{ij}^{m^*}b_{kt}^{m^*}(b_{tj}^{l^*}b_{ik}^{l^*}-b_{tk}^{l^*}b_{ij}^{l^*})\\
&+\sum_{ijkmtl}
b_{ij}^{m^*}b_{kt}^{m^*}\{b_{tj}^{l^*}c_{ik}^{l^*}+c_{tj}^{l^*}b_{ik}^{l^*}+c_{tj}^{l^*}c_{ik}^{l^*}
-b_{tk}^{l^*}c_{ij}^{l^*}-c_{tk}^{l^*}c_{ij}^{l^*}\}\\
&=c\|B\|^2+\sum_{ijkmtl}
b_{ij}^{m^*}b_{kt}^{m^*}(b_{tj}^{l^*}b_{ik}^{l^*}-b_{tk}^{l^*}b_{ij}^{l^*})
+\frac{2n}{n+2}\sum_{jklmt}b_{jk}^{m^*}b_{kt}^{m^*}b_{tj}^{l^*}H^{l^*}\\
&+ \sum_{ijkmtl}
b_{ij}^{m^*}b_{kt}^{m^*}(c_{tj}^{l^*}c_{ik}^{l^*}-c_{tk}^{l^*}c_{ij}^{l^*})\\
&=c\|B\|^2+\frac{n^2}{(n+2)^2}\|B\|^2|H|^2+\sum_{ijkmtl}
b_{ij}^{m^*}b_{kt}^{m^*}(b_{tj}^{l^*}b_{ik}^{l^*}-b_{tk}^{l^*}b_{ij}^{l^*})\\
&+ \frac{2n}{n+2}\sum_{jklmt}
b_{jk}^{m^*}b_{kt}^{m^*}b_{tj}^{l^*}H^{l^*}+\frac{2n^2}{(n+2)^2}\sum_{ijkm}
b_{ij}^{m^*}b_{jk}^{m^*}H^{i^*}H^{k^*}
\endaligned
$$
$$
\aligned
(II)&=(n-1)c\|B\|^2+\sum_{ijklmt}b_{ij}^{m^*}b_{ti}^{m^*}\{nb_{tj}^{l^*}H^{l^*}+nc_{tj}^{l^*}H^{l^*}
-b_{tk}^{l^*}b_{kj}^{l^*}-b_{tk}^{l^*}c_{kj}^{l^*}\\
& - c_{tk}^{l^*}b_{kj}^{l^*}-c_{tk}^{l^*}c_{kj}^{l^*}\}\\
&=(n-1)c\|B\|^2-\sum_{ijklmt}b_{ij}^{m^*}b_{ti}^{m^*}b_{tk}^{l^*}b_{kj}^{l^*}
+\frac{n^2-2n}{n+2}\sum_{ijlmt}
b_{ij}^{m^*}b_{ti}^{m^*}b_{tj}^{l^*}H^{l^*}\\
&+\frac{n^2}{n+2}\{\|B\|^2|H|^2+2\sum_{ijmt}b_{ij}^{m^*}b_{ti}^{m^*}H^{j^*}H^{t^*}\}
-\sum_{ijklmt}b_{ij}^{m^*}b_{ti}^{m^*}c_{tk}^{l^*}c_{kj}^{l^*}\\
&=(n-1)c\|B\|^2+\frac{n^3}{(n+2)^2}\|B\|^2|H|^2-\sum_{ijkmtl}
b_{ij}^{m^*}b_{ti}^{m^*}b_{tk}^{l^*}b_{kj}^{l^*}\\
&+ \frac{n^2-2n}{n+2}\sum_{ijlmt}
b_{ij}^{m^*}b_{ti}^{m^*}b_{tj}^{l^*}H^{l^*}+\frac{n^2(n-2)}{(n+2)^2}\sum_{ijmt}
b_{ij}^{m^*}b_{ti}^{m^*}H^{j^*}H^{t^*}\endaligned
$$
$$
\aligned
(III)&=\sum_{ijklmt}b_{ij}^{m^*}b_{ki}^{t^*}\{b_{lj}^{m^*}b_{lk}^{t^*}-b_{lk}^{m^*}b_{lj}^{t^*}
+b_{lj}^{m^*}c_{lk}^{t^*}+c_{lj}^{m^*}b_{lk}^{t^*}-b_{lk}^{m^*}c_{lj}^{t^*}-c_{lk}^{m^*}b_{lj}^{t^*}\\
& + c_{lj}^{m^*}c_{lk}^{t^*}-c_{lk}^{m^*}c_{lj}^{t^*}\}-c\|B\|^2\\
&= -c\|B\|^2-\frac{n^2}{(n+2)^2}\|B\|^2|H|^2+\sum_{ijklmt}
b_{ij}^{m^*}b_{ki}^{t^*}(b_{lj}^{m^*}b_{lk}^{t^*}-b_{lk}^{m^*}b_{lj}^{t^*})\\
& - \frac{2n^2}{(n+2)^2}\sum_{ijkm}
b_{ij}^{m^*}b_{jk}^{m^*}H^{i^*}H^{k^*}-\frac{2n}{n+2}\sum_{ijlmt}
b_{ij}^{m^*}b_{lj}^{m^*}b_{li}^{t^*}H^{t^*}\endaligned
$$
Set $B_{m^*}=(b_{ij}^{m^*})$. Then we get
$$
\aligned
\frac{1}{2}\triangle\|B\|^2&=\sum_{ijkm}{(b_{ijk}^{m^*})}^2+\sum_{ijkm}
b_{ij}^{m^*}b_{ijkk}^{m^*}\\
&=\sum_{ijkm}{(b_{ijk}^{m^*})}^2+(n+1)c\|B\|^2+\frac{n^2}{(n+2)}\|B\|^2|H|^2\\
&+\sum_{ijkmtl}
b_{ij}^{m^*}b_{kt}^{m^*}(b_{tj}^{l^*}b_{ik}^{l^*}-b_{tk}^{l^*}b_{ij}^{l^*})
-\sum_{ijkmtl}
b_{ij}^{m^*}b_{ki}^{t^*}(b_{lj}^{m^*}b_{lk}^{t^*}-b_{lk}^{m^*}b_{lj}^{t^*})\\
&-\sum_{ijkmtl}
b_{ij}^{m^*}b_{it}^{m^*}b_{tk}^{l^*}b_{kj}^{l^*}+n\sum_{ijlmt}
b_{ij}^{m^*}b_{jl}^{m^*}b_{li}^{t^*}H^{t^*}
+\frac{n^2}{n+2}\sum_{ijkm} b_{ij}^{m^*}b_{jk}^{m^*}H^{i^*}H^{k^*} \\
&=\sum_{ijkm}{(b_{ijk}^{m^*})}^2+(n+1)c\|B\|^2+\frac{n^2}{(n+2)}\|B\|^2|H|^2
+\sum_{ij} tr{(B_{i^*}B_{j^*}-B_{j^*}B_{i^*})}^2\\
&-\sum_{ij} (trB_{i^*}B_{j^*})^2+\frac{n^2}{n+2}\sum_{ijkm}
b_{ij}^{m^*}b_{jk}^{m^*}H^{i^*}H^{k^*}+n\sum_{ijmt}
b_{ij}^{m^*}b_{jl}^{m^*}b_{li}^{t^*}H^{t^*}\\
&\geq\sum_{ijkm}{(b_{ijk}^{m^*})}^2+(n+1)c\|B\|^2+\frac{n^2}{(n+2)}\|B\|^2|H|^2
-\sum_{i}S_{i^*}-2\sum_{i\neq j}S_{i^*}S_{j^*}\\
&+n|H|\sum_{i}S_{i^*}\lambda_{i}+\frac{n^2}{n+2}S_{1^*}|H|\endaligned
$$
where
$$
S_{i^*}=\sum_{jt} (b_{jt}^{i^*})^2,\quad b_{ij}^{1^*}=\lambda_{i}\delta_{ij},\quad S_{H}=\sum_{i}\lambda_{i}^2
$$
Here we assume $e_{1^*}//H$ and $S_H=\sum(b_{ij}^{1^*})^2$.

For a matrix $A=(a_{ij})$, we denote by $N(A)$ the square norm of
$A$ as in [3], i.e., $$N(A)=trace(A^t A)=\sum_{ij}(a_{ij})^2$$
\proclaim{Lemma 3.2}([6]) Let $A_{1},\cdots,A_{p}$ be symmetric
$(n\times n)$-matrices $(p\geq 2)$. Denote
$S_{\alpha\beta}=trace(A_{\alpha}A_{\beta})$,
$S_{\alpha}=S_{\alpha\alpha}=N(A_{\alpha})$, $S=\sum_{i=1}^p
S_{i}$. Then
$$\sum_{\alpha\beta}N(A_{\alpha}A_{\beta}-A_{\beta}A_{\alpha})+\sum_{\alpha\beta}S_{\alpha\beta}^2
\leq \frac{3}{2}S^2$$
\endproclaim

From Lemma 3.2, we know
$$
\aligned
\frac{1}{2}\triangle\|B\|^2&\geq\sum_{ijkm}{(b_{ijk}^{m^*})}^2+(n+1)c\|B\|^2+\frac{n^2}{(n+2)}\|B\|^2|H|^2
-\frac{n+2}{2}\sum_{i}S^2_{i^*}\\
&-2\sum_{i\neq j}S_{i^*} S_{j^*}+\frac{n}{2}\sum_{i}(\lambda_i^2|H|^2+2|H|S_{i^*}\lambda_i+S_{i^*}^2)\\
&\geq
(n+1)c\|B\|^2+\frac{n^2}{(n+2)}\|B\|^2|H|^2-\frac{n+2}{2}(\sum_{i}S_{i^*}^2+2\sum_{i\neq
j}S_{i^*}S_{j^*})\\
& +n\sum_{i\neq
j}S_{i^*}S_{j^*}+\frac{n}{2}\sum_{i}(\lambda_i|H|+S_{i^*})^2\\
&\geq(n+1)c\|B\|^2+\frac{n^2}{(n+2)}\|B\|^2|H|^2-\frac{3(n+2)}{4}\|B\|^4\endaligned
$$

\proclaim{Theorem 3.3} Let $\psi: M\longrightarrow
\widetilde{M}^n(4c)$ $(n\geq 2)$ be a (non-minimal) Lagrangian
immersion of an $n$-dimensional closed manifold $M$ into the
complex space form $\widetilde{M}^n(4c)$ $ (c \geq 0)$ with
conformal Maslov form. If
$$\|B\|^2\leq\frac{4(n+1)c}{3(n+2)}+\frac{4n^2|H|^2}{3(n+2)^2} $$
or equivalently,
$$\|h\|^2\leq\frac{4(n+1)c}{3(n+2)}+\frac{n^2(9n+22)|H|^2}{3(n+2)^2}$$ then
$\psi(M)$ is the Whitney sphere.
\endproclaim
\demo{Proof} Since $M$ is closed, we have
$$
\aligned 0=\int_M\frac{1}{2}\triangle\|B\|^2&\geq \int_M
\{(n+1)c\|B\|^2+\frac{n^2}{(n+2)}\|B\|^2|H|^2-\frac{3(n+2)}{4}\|B\|^4\}\\
&=\int_M
\|B\|^2\{(n+1)c+\frac{n^2}{(n+2)}|H|^2-\frac{3(n+2)}{4}\|B\|^2\}\endaligned
$$
By the assumption
$$\|B\|^2\leq\frac{4(n+1)c}{3(n+2)}+\frac{4n^2|H|^2}{3(n+2)^2}$$
thus we get $B=0$ on $M$, i.e.,
$$h(X,Y)=\frac{n}{n+2}\{g(X,Y)H+g(JX,H)JY+g(JY,H)JX\}$$
for any tangent vector fields $X$, $Y$ of $M$. It follows from
Lemma 2.3 that $\psi(M)$ is the  Whitney sphere. From (3), we see
that the condition for $\|B\|^2$ is equivalent to
$$\|h\|^2\leq\frac{4(n+1)c}{3(n+2)}+\frac{n^2(9n+22)|H|^2}{3(n+2)^2}$$
This complete the proof of the Theorem.\qed
\enddemo
\remark{Remark 3.1}(a) The non-minimality condition is void for
the case $c=0$ . Obviously (3) tells us that ([7],[2])
$$\|h\|^2\geq \frac{3n^2|H|^2}{n+2}$$

(b) Let $M=S^1(r_1)\times\dots\times S^1(r_n)\hookrightarrow
C^1\times\dots\times C^1=C^n$ be the standard Lagrangian flat
torus in $C^n$, which has parallel mean curvature vector. By a
direct computation, we know that it has conformal Maslov form and
$$\|B\|^2=\frac{n^2(n-1)H^2}{n+2}$$ We also have the
standard flat Lagrangian torus in $\widetilde{M}^n(4c)$ with $c>0$
(see for example \cite{Ch}), which has parallel mean curvature,
and thus has conformal Maslov form. These examples show that some
upper bound for $\|B\|^2$ is necessary to characterize Whitney
spheres. It would be interesting to get the optimal gap theorem.
\endremark

\proclaim{Corollary 3.4} Let $\psi: M\longrightarrow
\widetilde{M}^n(4c)$ $(n\geq 2)$ be a (non-minimal) Lagrangian
immersion of an $n$-dimensional closed manifold $M$ into complex
space form $\widetilde{M}^n(4c)$ $ (c > 0)$ with conformal Maslov
form. If
$$\|B\|^2\leq\frac{4(n+1)c}{3(n+2)}$$  then
$\psi(M)$ is the Whitney sphere.
\endproclaim

\vskip 1.0 true cm 

\vskip 3.8 true cm

\noindent Xiaoli Chao

\noindent Department of Mathematics,

\noindent Southeast University, Nanjing, 210096,

\noindent P. R. China

\noindent Email address: xlchao\@seu.edu.cn \vskip 0.2 true cm

\noindent and \vskip 0.2 true cm

\noindent Yuxin Dong\noindent

\noindent Institute of Mathematics,

\noindent Fudan University, Shanghai, 200433,

\noindent P. R. China

\noindent And

\noindent Key Laboratory of Mathematics\noindent for Nonlinear
Sciences (Fudan University),

\noindent Ministry Education

\noindent Email address: yxdong\@fudan.edu.cn

\Refs \widestnumber\key{SW1}

\ref\key 1 \by B. Y. Chen\paper Riemannian geometry of Lagrangian
submanifolds\paperinfo Taiwanese J. of  Math., {\bf{5}}(2001),
681-723
\endref

\ref\key 2\by B. Y. Chen\paper Jacobi's elliptic functions and
Lagrangian immersions\paperinfo Proc. Royal Soc. Edin., {\bf{126}}
(1996), 687-704\endref

\ref\key 3\by S. S. Chern, M. do Carmo and S. Kobayashi\paper
Minimal submanifolds of a sphere with second fundamental form of
constant length\paperinfo Functional Analysis and Related Fields (
F. E. Browder ed.), Springer-Verlag, New York  (1970)\endref

\ref\key 4\by I.Castro, C.R.Montealegre and F.Urbano\paper Closed
conformal vector fields and Lagrangian submanifolds in complex
space forms\paperinfo Pacific J. Math., {\bf{199}} (2001),
269-301\endref

\ref\key 5\by M. L. Gromov\paper Psedoholomorphic curves in
symplectic manifolds\paperinfo invent. Math., {\bf{82}} (1985),
307-347\endref

\ref\key 6\by A.M.Li and J.M.Li\paper An intrinsic theorem for
minimal submanifolds in a sphere\paperinfo Archiv der Math.,
{\bf{58}} (1992), 582-594\endref

\ref\key 7\by  A. Ros and F.Urbano\paper Lagrangian submanifolds
of $C^n$ with conformal Maslov form and the Whitney
sphere\paperinfo J. Math. Soc. Japan , {\bf{50}} (1998),
203-226\endref

\ref\key 8\by A. Weinstein\paper Lectures on symplectic
manifolds\paperinfo Conference board of the Mathematical
Scientific, {\bf{29}}(1977)\endref

\endRefs
\enddocument